\documentclass[12pt]{amsart}

\usepackage{latexsym, graphics, psfrag, amscd, amssymb, graphicx, psfrag}
\newtheorem{thm}{Theorem}
\newtheorem{tm}[thm]{Theorem}

\newtheorem{lm}[thm]{Lemma}
\newtheorem{rem}[thm]{Remark}
\newtheorem{cor}[thm]{Corollary}
\newtheorem{pro}[thm]{Proposition}
\newtheorem{de}[thm]{Definition}

\newcommand{\ra}{\rightarrow}
\newcommand{\Si}{\Sigma}
\newcommand{\tS}{\tilde{\Sigma}}
\newcommand{\M}{\mathcal{M}}
\newcommand{\Z}{\mathbb{Z}}

\newcommand{\A}{\mathcal A}
\newcommand{\cG}{\mathcal G}

\newcommand{\g}{\mathfrak g}

\begin{document}

\parskip=0.35\baselineskip
\baselineskip=1.2\baselineskip
\title[Nonorientable surface and its double cover]{The real locus
of an involution map on the moduli space of flat connections on a
Riemann surface}

\author{Nan-Kuo Ho}
\address{Department of Mathematics\\ National Cheng-Kung University}\thanks{The author
  was supported by grant NSC 92-2119-M-006-006 and grant from OGS and OGSST}
\email{nankuo@math.toronto.edu\\ nankuo@mail.ncku.edu.tw}


\subjclass{53}
\date{\today}

\begin{abstract}
It is known that every nonorientable surface $\Sigma$ has an
orientable double cover $\tilde{\Sigma}$. The covering map induces
an involution on the moduli space $\tilde{\M}$ of gauge
equivalence classes of flat $G$-connections on $\tilde{\Sigma}$.
We identify the relation between the moduli space $\M$ and the
fixed point set of the moduli space $\tilde{\M}$. In particular,
$\M$ is isomorphic to the fixed point set of $\tilde{\M}$ if and
only if the order of the center of $G$ is odd. One important
application is that we give a way to construct a minimal
Lagrangian submanifold of the moduli space $\tilde{\M}$.

\end{abstract}

\maketitle

\section{Introduction}

We know that every nonorientable surface $\Sigma$ has an
orientable double cover $\tilde{\Sigma}$. The covering map induces
an involution on the space of all flat $G$-connections on
$\tilde{\Sigma}$ such that the space of all flat $G$-connections
on $\Sigma$ sits in it as the fixed point set of the induced
involution. Again, there is an induced involution on the moduli
space $\tilde{\M}$ of gauge equivalence classes of flat
$G$-connections on $\tilde{\Sigma}$. However, we do not know if we
will still have that the moduli space $\M$ of flat $G$-connections
on $\Si$ sits in $\tilde{\M}$ as the fixed point set of the
involution since now we take the quotient by gauge transformations
on both spaces.

We show in this paper that when the order of the center of $G$ is odd,
there exists an
isomorphism between (the smooth part which is an open dense set of) the
moduli space $\M$ and (the smooth part which is an open dense set of)
the fixed point set of the involution on the
moduli space $\tilde{\M}$.

The benefit of this result is that since
the moduli space of gauge equivalence classes of flat
$G$-connections on an orientable surface is K\"{a}hler, the fixed
point set of this involution, which is an anti-symplectic,
anti-holomorphic isometry, is a totally geodesic Lagrangian
submanifold of the K\"{a}hler manifold.

We constructed a model similar to the standard model used to
identify the moduli space of flat $G$-connections on a surface
with the space of representations from the fundamental group of
the surface into $G$, divided by the conjugate action (for
example:\cite{G2}) by introducing two base points on the surface
instead of one base point. This model is constructed in such a way
that it is clear how exactly the involution works on both the
surface $\tilde{\Si}$ and the moduli space $\tilde{\M}$.

Through out this paper, denote $\Si$ a compact, closed surface;
$G$ a compact, connected, semisimple Lie group; and $e$ the
identity element of the Lie group $G$. Our main result is, (note
that the smooth part of the moduli space is an open dense subset)

\begin{tm}
Let $\Si$ be a compact, closed nonorientable surface and
$\tilde{\Si}$ its orientable double cover. Then the natural map from
the moduli space of flat connections on $\Si$, $\M$, into
the fixed point set of the involution on the
moduli space of flat connections on $\tilde{\Si}$, $\tilde{\M}$, is
a $\mid Z(G)/2Z(G)\mid$ to one map. In particular,
(the smooth part of) the moduli space $\M$ is
isomorphic to (the smooth part of) the fixed point set of the
involution on the moduli space $\tilde{\M}$
if and only if $\mid Z(G)\mid$ is odd.
\end{tm}

The paper is organized as follows: section two explains what is
the fixed point set of the involution on $\tilde{\M}$; section
three gives the relation between the fixed point set of
$\tilde{\M}$ and $\M$; Examples are given in section four, and the
last section gives an application of the theorem to the moduli
space of semistable vector bundles on a Riemann surface.

\textit{Acknowledgment:} This article is based on part of the
author's Ph.D. thesis under the supervision of Professor Lisa
Jeffrey. It is a pleasure to thank Chiu-Chu Melissa Liu, and
Professor Eckhard Meinrenken for sharing some of their expertise
over hours of conversations and their impact on this article. The
author also wants to thank Professor Chris Woodward for valuable
suggestions for last section.

\section{The fixed point set of the involution of $\tilde{\M}$}

Consider $\Si$, a nonorientable surface, and $\tS$ its orientable
double cover. Let $\M$ denote the moduli space of flat
$G$-connections on $\Si$ and $\tilde{\M}$ denote the moduli space
of flat $G$-connections on $\tS$. The $\Z_2$ action on $\tS$ gives
the quotient space $\Si$, and this $\Z_2$-action induces an
involution $\tau$ on the space of all $\g$-valued one forms on
$\tS$ and the fixed point set of the space of of all $\g$-valued
one forms on $\tS$ is equivalent to the space of all $\g$-valued
one forms on $\Si$.

Let $\tau: \Omega^1(\tS,\g) \ra \Omega^1(\tS,\g)$ denote the
involution ($\Z_2$-action) and $A$ is a one form on $\tilde{\Si}$.
Then $A$ is a fixed point of $\tau$ if $\tau(A)(x)=A(x)$ for all
$x \in \tS$, i.e. $A(\tau(x))=A(x)$, i.e. $A$ has the same value
at $x$ and $\tau(x)$ which means that $A$ corresponds to a one
form in $\Omega^1(\Si,\g)$.

Thus, the the fixed point set of the involution on the flat
$G$-connections on the orientable double cover is the same as the
space of all flat $G$-connections on the nonorientable surface,
i.e. $(\mathcal{A}_{flat}(\tS))^{\tau}=\mathcal{A}_{flat}(\Si)$.

Fix a metric $h$ on $\Si$ and let $\pi: \tS \ra \Si$ be the
covering map. Let $\pi^{\ast}h$ be the (corresponding) pullback
metric on $\tS$ then $\pi^{\ast}h$ determines a complex structure
$J$ on $\tilde{\M}$ such that $(\tilde{\M},\omega)$ is a
K\"{a}hler manifold. In our situation here the $J$ is the Hodge
star $\ast$ on $\Omega^1(\tS,\g)$, and
 $\omega$ is the usual symplectic structure on $\tilde{\M}$
i.e. $\omega(a,b)=\int_{\tS} Tr(a\wedge b)$ where $Tr$ stands for
the Killing form on the Lie algebra $\g$.

We have the following two observations:

\begin{lm}
The isometric involution $\tau:\A_{flat}(\tS) \ra \A_{flat}(\tS)$
is anti-holomorphic.
\end{lm}\paragraph{Proof}

Since $\tau$ is an orientation reversing isometry, we have
$<a,b>=<\tau_{\ast}a,\tau_{\ast}b>$ and $a\wedge
b=-\tau_{\ast}a\wedge \tau_{\ast}b$. Thus
\begin{eqnarray*}
<a,b>&=:&\int Tr(a\wedge Jb)=
\int -Tr(\tau_{\ast}a\wedge \tau_{\ast}Jb) \\
<\tau_{\ast}a,\tau_{\ast}b>&=:& \int Tr(\tau_{\ast}a\wedge
J\tau_{\ast}b)
\end{eqnarray*}
for all $a,b \in T_{A}\A_{flat}(\tS)$, thus we have
$J\tau_{\ast}=-\tau_{\ast}J$ which means $\tau$ is
anti-holomorphic. $\Box$

\begin{pro}
Let $M$ be a K\"{a}hler manifold. If $f:M \ra M$ is an
anti-holomorphic isometric involution on $M$, and if the fixed
point set of $f$, $N$, is nonempty, then $N$ is a totally
geodesic, totally real, Lagrangian submanifold.
\end{pro}

\paragraph{Proof}
We refer to ~\cite{Oh}.

(1) We would like to show that a fixed point set of an isometry is
a totally geodesic submanifold.

If $N$ is the fixed point set of an isometry $f:M \ra M$, and
suppose $r(t)$ is a geodesic in $M$ where $r(0)$ is in $N$ and
$r'(0) \in T_{r(0)}N$, then $f\circ r(0)=r(0)$ since $r(0) \in N$
and $(f\circ r)'(0)=f_{\ast}r'(0)=r'(0)$ since $r'(0) \in
T_{r(0)}N$. Thus $f\circ r$ is a geodesic (since f is an isometry
and $r$ is a geodesic) with $f\circ r(0)=r(0)$ and $f\circ
r'(0)=r'(0)$. It follows that $r(t)=f\circ r(t)$ by the uniqueness
of geodesics, which means $r(t)$ is fixed by $f$ so $r(t) \in N$.
This shows that $N$ is a totally geodesic submanifold.

(2) The fixed point set of an anti-holomorphic involution is
totally real.

Suppose $f$ is the anti-holomorphic involution, and $N$ is the
fixed point set. Let $M$ and $J$ be given. Now since $f$ is
anti-holomorphic, $Jf_{\ast}=-f_{\ast}J$. Now, if $v \in T_xN$
then $f_{\ast}v=v$ and $Jv=Jf_{\ast}v=-f_{\ast}(Jv)$ so $Jv$ is
not fixed by $f_{\ast}$. Thus $Jv \notin T_xN$, so $T_xN\cap
J(T_xN)=0$ and $N$ is totally real. $\Box$

Since the moduli space of flat $G$-connections on a Riemann
surface is a K\"{a}hler manifold (ref: ~\cite{AB}), and
$\tilde{\M}^{\tau}$ is nonempty ( since the set of
$(\mathcal{A}_{flat}(\tS))^{\tau}$ is nonempty) we have the
following corollary:

\begin{cor}
The real locus of the involution on the moduli space $\tilde{\M}$
is a totally real, totally geodesic Lagrangian submanifold of
$\tilde{\M}$.
\end{cor}

\section{The map between the fixed point set of the involution of $\tilde{\M}$ and $\M$} \label{projective}

In this section we discuss the relation between $\M$ and the fixed
point set of $\tilde{\M}$ according to its topological type of the
base space $\Si$.

{\bf Case 1. $\Si= \Si_\ell \sharp RP^2$}

Let $\Si$ be a compact nonorientable surface whose topological
type is the connected sum of a Riemann surface of genus $\ell$ and
the real projective plane. Let $\tS$ be the orientable double
cover of $\Si$. Then $\tS$ is a compact orientable surface of
genus $2\ell$.

Let us look at the picture below,

\begin{figure}[!htbp]
        \begin{center}
                \psfrag{a2}[c][c][1][0]{$\bar{A}$}
                \psfrag{b2}[c][c][1][0]{$\bar{B}$}
                \psfrag{c2}[c][c][1][0]{$\bar{C}$}
                \psfrag{c1}[c][c][1][0]{$C$}
                \psfrag{a1}[c][c][1][0]{$A$}
                \psfrag{b1}[c][c][1][0]{$B$}
                \psfrag{p1}[c][c][1][0]{$P_{+}$}
                \psfrag{p2}[c][c][1][0]{$P_{-}$}
                \includegraphics[scale = 0.7]
                {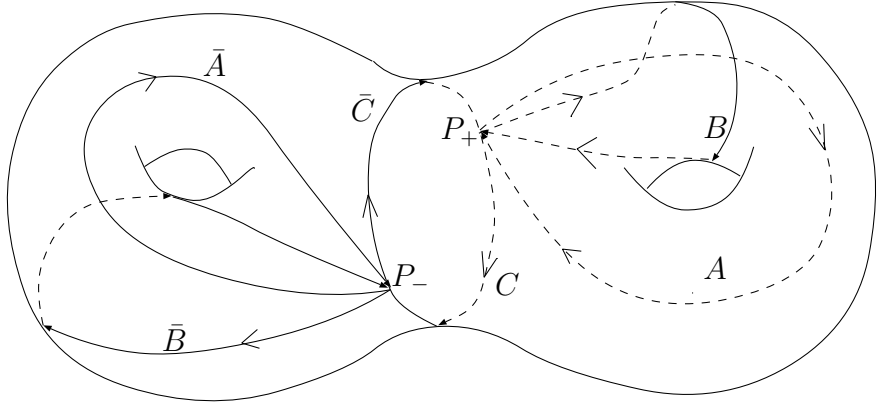}
        \end{center}
        \caption{Holonomy on the double cover}
\end{figure}

where $a_i,b_i$ (resp. $\bar{a_i},\bar{b_i}$) are the holonomies
of flat connections around loops $A_i$ and $B_i$ (resp.
$\bar{A_i}$ and $\bar{B_i}$) based at $P_+$ (resp. $P_-$), while
$c$ is the parallel transport along an arc from $P_+$ to $P_-$ and
$\bar{c}$ is the parallel transport along an arc from $P_-$ to
$P_+$.

We can then define the moduli space of flat connections on
$\tilde{\Si}$ by introducing two base points $P_{+},P_{-}$,
\begin{eqnarray*}
\M(\tilde{\Si})&=&[\A_{F}(\tilde{\Si})/\cG(\tilde{\Si},P_{+},P_{-})]/G^2\\
&=&\left\{(a_1,b_1,\cdots,a_\ell,b_\ell,c,
    \bar{a}_1,\bar{b}_1,\cdots,\bar{a}_\ell,\bar{b}_\ell,\bar{c})
\in G^{2(2\ell+1)}\mid \right.\\
&& \left. \prod_{i=1}^{\ell}[a_i,b_i]=c\bar{c}
,~\prod_{i=1}^\ell[\bar{a}_i,\bar{b}_i]=\bar{c}c \right\}/G\times
G
\end{eqnarray*}
where $G\times G$ acts as \begin{eqnarray*}
&& (g_1,g_2)\cdot(a_1,b_1,\cdots,a_\ell,b_\ell,c,
\bar{a}_1,\bar{b}_1,\cdots,\bar{a}_\ell,\bar{b}_\ell,\bar{c})
\\&=&
(g_1a_1g_1^{-1},g_1b_1g_1^{-1},\cdots,g_1a_\ell g_1^{-1},g_1b_\ell
g_1^{-1}, g_1 c g_2^{-1},\\&&
 g_2\bar{a}_1g_2^{-1},g_2\bar{b}_1g_2^{-1},\cdots,
  g_2\bar{a}_\ell g_2^{-1},g_2\bar{b}_\ell g_2^{-1},g_2\bar{c}g_1^{-1}).
\end{eqnarray*}

and the moduli space of flat connections on $\Si$,
\begin{eqnarray*}
\M(\Si)&=&[\A_{F}(\Si)/\cG(\Si,P)]/G\\
&=&\left\{(a_1,b_1,\cdots,a_\ell,b_\ell,c) \in G^{2\ell+1} |
\prod_{i=1}^{\ell}[a_i,b_i]=c^2 \right\}/G
\end{eqnarray*}
where $G$ action is the conjugate action.

{\bf Case 2. $\Si=\Si_\ell\sharp ~\mbox{Klein bottle}$.}

Let $\Si$ be a compact nonorientable surface whose topological
type is the connected sum of a Riemann surface of genus $\ell$ and
a Klein bottle. Let $\tS$ be the orientable double cover of $\Si$.
Then $\tS$ is a compact orientable surface of genus $2\ell+1$.

\begin{figure}[!htbp]
        \begin{center}
                \psfrag{a2}[c][c][1][0]{$\bar{A}$}
                \psfrag{b2}[c][c][1][0]{$\bar{B}$}
                \psfrag{a1}[c][c][1][0]{$A$}
                \psfrag{b1}[c][c][1][0]{$B$}
                \psfrag{p2}[c][c][1][0]{$P_{-}$}
                \psfrag{p1}[c][c][1][0]{$P_{+}$}
                \psfrag{x2}[c][c][1][0]{$\bar{C}$}
                \psfrag{y2}[c][c][1][0]{$\bar{D}$}
                \psfrag{x1}[c][c][1][0]{$C$}
                \psfrag{y1}[c][c][1][0]{$D$}
                \includegraphics[scale = 0.7]
                {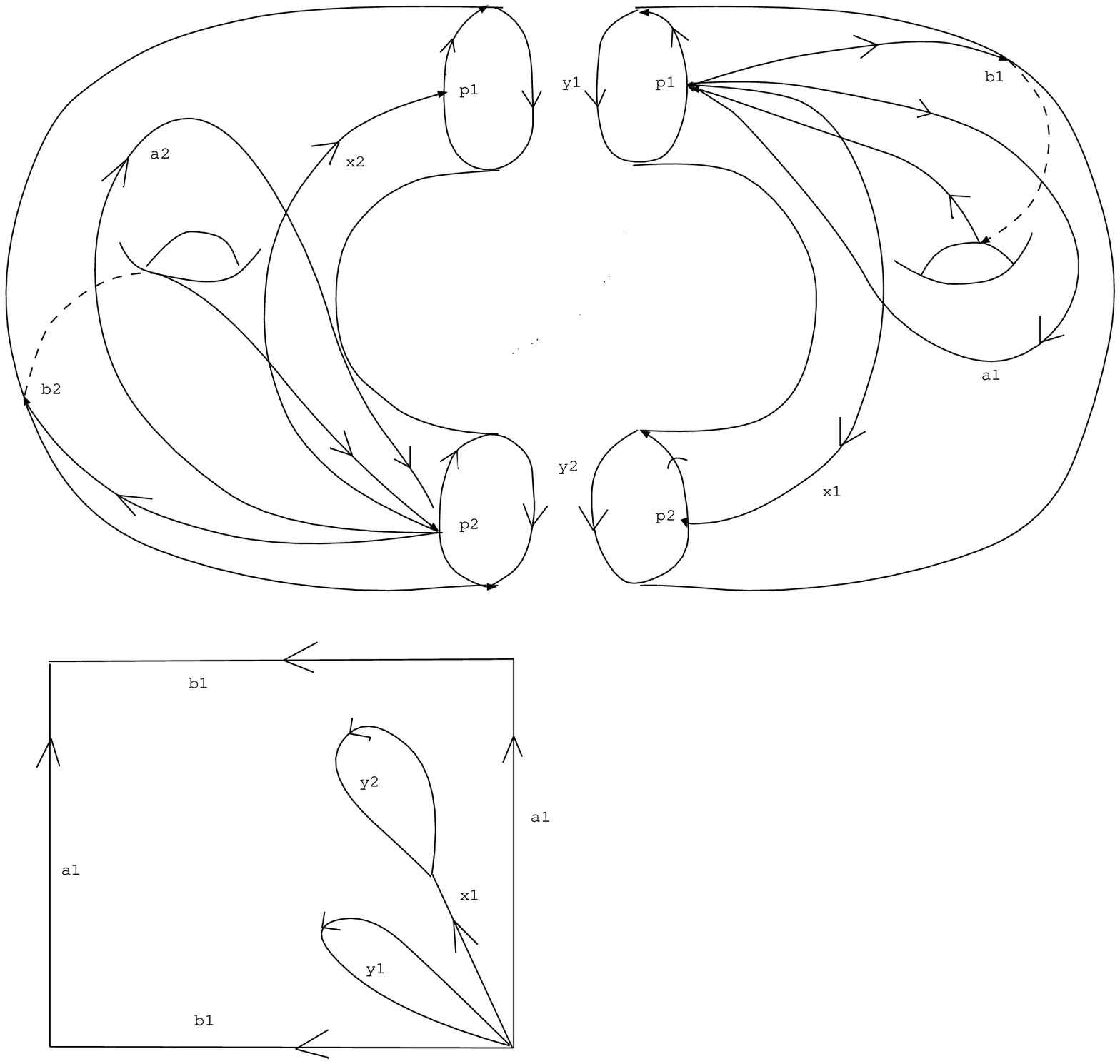}
        \end{center}
        \caption{Holonomy on the double cover $\Si_3$}
\end{figure}
Figure 2 is the picture of holonomy on the double cover $\Si_3$,
where $a_i,b_i,d$ (resp. $\bar{a_i},\bar{b_i},\bar{d}$) are the
holonomies of flat connections around loops $A_i$, $B_i$, and $D$
(resp. $\bar{A_i}$, $\bar{B_i}$, and $\bar{D}$) based at $P_+$
(resp. $P_-$), while $c$ is the parallel transport along an arc
from $P_+$ to $P_-$ and $\bar{c}$ is the parallel transport along
an arc from $P_-$ to $P_+$.

This means
\[aba^{-1}b^{-1}d^{-1}c\bar{d}^{-1}c^{-1}=\mbox{id.}\]thus
\[aba^{-1}b^{-1}=c\bar{d}c^{-1}d\]

We can then define the moduli space of flat connections on
$\tilde{\Si}$ by introducing two base points $P_{+},P_{-}$
\begin{eqnarray*}
\M(\tilde{\Si})&=&[\A_{F}(\tilde{\Si})/\cG(\tilde{\Si},P_{+},P_{-})]/G^2\\
&=&\left\{(a_1,b_1,\cdots,a_\ell,b_\ell,d,c,
    \bar{a}_1,\bar{b}_1,\cdots,\bar{a}_\ell,\bar{b}_\ell,\bar{d},\bar{c})
\in G^{2(2\ell+2)}\mid \right.\\
&&\ \left. \prod_{i=1}^{\ell}[a_i,b_i]=c\bar{d}c^{-1}d
,\prod_{i=1}^\ell[\bar{a}_i,\bar{b}_i]=\bar{c}d\bar{c}^{-1}\bar{d}
  \right\}/G \times G
\end{eqnarray*}where $G\times G$ acts as
\begin{eqnarray*}
&& (g_1,g_2)\cdot(a_1,b_1,\cdots,a_\ell,b_\ell,d,c,
\bar{a}_1,\bar{b}_1,\cdots,\bar{a}_\ell,\bar{b}_\ell,\bar{d},\bar{c})\\
&=& (g_1a_1g_1^{-1},g_1b_1g_1^{-1},\cdots,g_1a_\ell
g_1^{-1},g_1b_\ell g_1^{-1}, g_1dg_1^{-1},g_1 c g_2^{-1},\\&&
g_2\bar{a}_1g_2^{-1},g_2\bar{b}_1g_2^{-1},\cdots,
  g_2\bar{a}_\ell g_2^{-1},g_2\bar{b}_\ell g_2^{-1},g_2\bar{d}g_2^{-1},g_2\bar{c}g_1^{-1}).
\end{eqnarray*}

and the moduli space of flat connections on $\Si$,
\begin{eqnarray*}
\M(\Si)&=&[\A_{F}(\Si)/\cG(\Si,P)]/G\\
&=&\left\{(a_1,b_1,\cdots,a_\ell,b_\ell,d,c) \in G^{2\ell+2} |
\prod_{i=1}^{\ell}[a_i,b_i]=cdc^{-1}d \right\}/G
\end{eqnarray*}where $G$ action is the conjugate action.

Thus the natural involution from $\M(\Si)$ to $\M(\tS)$ is either
\begin{eqnarray*}
\M(\Si)&\ra& \M(\tS)\\ ~[(a_1,b_1,\cdots,c)] &\mapsto &
[(a_1,b_1,\cdots,c,a_1,b_1,\cdots,c)]
\end{eqnarray*}
or
\begin{eqnarray*}
\M(\Si)&\ra& \M(\tS)\\ ~[(a_1,b_1,\cdots,d,c)] &\mapsto&
[(a_1,b_1,\cdots,d,c,a_1,b_1,\cdots,d,c)]
\end{eqnarray*}

We would like to introduce some notations for later convenience:

{\bf Case 1: $\Si= \Si_\ell \sharp RP^2$}

Set
\[
M=\left\{(a_1,b_1,\cdots,a_\ell,b_\ell,c,
    \bar{a}_1,\bar{b}_1,\cdots,\bar{a}_\ell,\bar{b}_\ell,\bar{c})
\in G^{2(2\ell+1)}\mid \right.\]\begin{equation}\label{M1}
 \left. \prod_{i=1}^{\ell}[a_i,b_i]=c\bar{c}
,~\prod_{i=1}^\ell[\bar{a}_i,\bar{b}_i]=\bar{c}c \right\}
\end{equation}

$K=G\times G$ acts on $M$ by
\[
 (g_1,g_2)\cdot(a_1,b_1,\cdots,a_\ell,b_\ell,c,
\bar{a}_1,\bar{b}_1,\cdots,\bar{a}_\ell,\bar{b}_\ell,\bar{c})=
(g_1a_1g_1^{-1},g_1b_1g_1^{-1},\cdots,\]
\begin{equation}\label{K1}g_1a_\ell
g_1^{-1},g_1b_\ell g_1^{-1}, g_1 c g_2^{-1},
 g_2\bar{a}_1g_2^{-1},g_2\bar{b}_1g_2^{-1},\cdots,
  g_2\bar{a}_\ell g_2^{-1},g_2\bar{b}_\ell g_2^{-1},g_2\bar{c}g_1^{-1})
\end{equation}

There is an involution $\tau:K\ra K$ given by
\begin{equation}\label{tK1}
\tau(g_{1},g_{2})=(g_{2},g_{1}),
\end{equation}
and an involution $\tau: M\ra M$ given by
\begin{equation}\label{tM1}
\tau(a_1,b_1,\cdots,a_\ell,b_\ell,c,
     \bar{a}_1,\bar{b}_1,\cdots,\bar{a}_\ell,\bar{b}_\ell,\bar{c})
    =(\bar{a}_{1},\bar{b}_{1},\cdots,\bar{a}_\ell,\bar{b}_\ell,\bar{c},
      a_1,b_1,\cdots,a_\ell,b_\ell,c).
\end{equation}

{\bf Case 2:  $\Si=\Si_\ell\sharp ~\mbox{Klein bottle}$}

Set \[ M=\left\{(a_1,b_1,\cdots,a_\ell,b_\ell,d,c,
    \bar{a}_1,\bar{b}_1,\cdots,\bar{a}_\ell,\bar{b}_\ell,\bar{d},\bar{c})
\in G^{2(2\ell+2)}\mid \right. \]\begin{equation}\label{M2}
 \left. \prod_{i=1}^{\ell}[a_i,b_i]=c\bar{d}c^{-1}d
,\prod_{i=1}^\ell[\bar{a}_i,\bar{b}_i]=\bar{c}d\bar{c}^{-1}\bar{d}
\right\}
\end{equation}

$K=G\times G$ acts on $M$ by
\begin{eqnarray*}
&& (g_1,g_2)\cdot(a_1,b_1,\cdots,a_\ell,b_\ell,d,c,
\bar{a}_1,\bar{b}_1,\cdots,\bar{a}_\ell,\bar{b}_\ell,\bar{d},\bar{c})\\
&=& (g_1a_1g_1^{-1},g_1b_1g_1^{-1},\cdots,g_1a_\ell
g_1^{-1},g_1b_\ell g_1^{-1}, g_1dg_1^{-1},g_1 c
g_2^{-1},\end{eqnarray*}\begin{equation}\label{K2}
g_2\bar{a}_1g_2^{-1},g_2\bar{b}_1g_2^{-1},\cdots,
  g_2\bar{a}_\ell g_2^{-1},g_2\bar{b}_\ell g_2^{-1},g_2\bar{d}g_2^{-1},g_2\bar{c}g_1^{-1}).
\end{equation}

There is an involution $\tau:K\ra K$ given by
\begin{equation}\label{tK2}
\tau(g_{1},g_{2})=(g_{2},g_{1}),
\end{equation}
and an involution $\tau: M\ra M$ given by
\[\tau(a_1,b_1,\cdots,a_\ell,b_\ell,d,c,
     \bar{a}_1,\bar{b}_1,\cdots,\bar{a}_\ell,\bar{b}_\ell,\bar{d},\bar{c})
\]\begin{equation}\label{tM2}
=(\bar{a}_{1},\bar{b}_{1},\cdots,\bar{a}_\ell,\bar{b}_\ell,\bar{d},\bar{c},
      a_1,b_1,\cdots,a_\ell,b_\ell,d,c).
\end{equation}

In both cases, we have $\tau(k\cdot x)=\tau(k)\cdot\tau(x)$ for
$k\in K$, $x\in M$ and it induces an involution $\tau: M/K\ra
M/K$,
\[
\tau([x])=[\tau(x)],
\]
which can be easily checked to be well-defined.

Notice that the action defined by equation (2) and (6) is just a
combination of conjugate action and an action as follows:
$(g,h)\cdot(x,y)=(gxh^{-1},hyg^{-1})$. Thus equation (2) and (6)
can be rewritten as
$$(g_1,g_2)\cdot(V,c,\bar{V},\bar{c})=(g_1Vg_1^{-1},g_1cg_2^{-1},g_2\bar{V}g_2^{-1},g_2\bar{c}g_1^{-1})$$
 if we denote $(a_1,b_1,\cdots,a_\ell,b_\ell)$ by $V$, and
$(ha_1h^{-1},hb_1h^{-1},\cdots,ha_\ell h^{-1},hb_\ell h^{-1})$ by
$hVh^{-1}$ for equation (2); $(a_1,b_1,\cdots,a_\ell,b_\ell,d)$ by
$V$, and $(ha_1h^{-1},hb_1h^{-1},\cdots,ha_\ell h^{-1},hb_\ell
h^{-1},hdh^{-1})$ by $hVh^{-1}$ for equation (6).

\begin{rem}
\end{rem}
With the definition of $M,K,\tau$ given in (\ref{M1}), (\ref{K1}),
(\ref{tK1}) and (\ref{tM1}), or (\ref{M2}), (\ref{K2}),
(\ref{tK2}), and (\ref{tM2}), the moduli space $\M(\tS)$ of gauge
equivalence classes of flat $G$-connections on $\tS$ is identified
with $M/K$, while the moduli space $\M(\Si)$ of gauge equivalence
classes of flat $G$-connections on $\Si$ is identified with
$M^\tau/K^\tau$. Also, we have the same action formulas for Case 1
and Case 2 if we use the simplified notation
$(V,c,\bar{V},\bar{c})$ to write points in $M$ where $V$
represents either $(a_1,b_1,\cdots,a_\ell,b_\ell)$ or
$(a_1,b_1,\cdots,a_\ell,b_\ell,d)$. {\bf Thus in this section, we
will be using $M,K,\tau$ and $(V,c,\bar{V},\bar{c})$ to show the
relation between $\M(\Si)$ and $\M(\tS)$.}


There is a natural map $I$ from the moduli space $M^\tau/K^\tau$
of flat $G$-connections on $\Si$ to $(M/K)^{\tau}$ the fixed point
set of the involution on the moduli space of flat $G$-connections
on $\tS$,
\begin{eqnarray*}
I:M^\tau/K^\tau &\ra& (M/K)^\tau \\
 ~[[(V,c,V,c)]] &\mapsto& [(V,c,V,c)]
\end{eqnarray*}

Here $[[.]]$ represents the equivalence class in $M^\tau/K^\tau$
and $[.]$ represents the equivalence class in $(M/K)^\tau$.

This map is well defined: If $[[(V,c,V,c)]]= [[(V',c',V',c')]]$ in
$M^\tau/K^\tau$ then there exists $(g,g)$ in $K^\tau$ such that
$(g,g)\cdot(V,c,V,c)=(V',c',V',c')$. Of course $[(V,c,V,c)] \in
(M/K)^\tau$ and $ [(V',c',V',c')] \in (M/K)^\tau$. Let $k=(g,g)$
thus $k \cdot(V,c,V,c)=(V',c',V',c')$ and they are the same
element in $(M/K)^\tau$.

We want to know if the map $I$ is surjective or injective.
\subsection{Surjectivity}\label{klein}

We know that points which have center as their stabilizer
correspond to the smooth part of the moduli space of flat
connections which is open dense( ref: \cite{G2} and \cite{AB}).
For this reason, we give the following definition:

\begin{de}\label{gen}
A point $x$ is called generic if its stabilizer is equal to the
center.
\end{de}

Suppose that $ x=(V,c,\bar{V},\bar{c})\in M$ satisfies $[x]\in
(M/K)^\tau$. Then there exists $k=(g_1,g_2)\in K$ such that
$k\cdot x=\tau(x)$. Recall that the action $k\cdot x$ is the
conjugate action on the $V, \bar{V}$ part and $
(g_1,g_2)\cdot(c,\bar{c})=(g_1cg_2^{-1},g_2\bar{c}g_1^{-1})$, and
$\tau(x)$ just switches $(V,c)$ and $(\bar{V},\bar{c})$. Thus we
have
\[
(g_1Vg_1^{-1},g_1cg_2^{-1},g_2\bar{V}g_2^{-1},g_2\bar{c}g_1^{-1})=(\bar{V},\bar{c},V,c)\]
$\Longleftrightarrow$
\[(g_2g_1Vg_1^{-1}g_2^{-1},g_2g_1cg_2^{-1}g_1^{-1})=(V,c)\]

For generic $[x]\in (M/K)^\tau$ we may assume that
$r^{-1}=g_2g_1\in Z(G)$. Then $k=(g,r^{-1}g^{-1})$, where
$g=g_1\in G$, $r\in Z(G)$ and we have $$x=(V,c,  gVg^{-1},
gcgr).$$ Let us define
\[
\tilde{x}:= (e,g^{-1})\cdot x =(V,cg,V,cgr).
\]
Note that $[\tilde{x}]=[x]$, and $\tilde{x}\in N_r$, for $N_r$ is
the subset of $M$ defined by

Case 1: $\Si= \Si_\ell \sharp RP^2$
\begin{eqnarray*}
N_r&=&\left\{(V,c,V,cr)=(a_1,b_1,\cdots,a_\ell,b_\ell,c,a_1,b_1,\cdots,a_\ell,b_\ell,cr)\mid\right.\\
&&\ \left. (a_1,b_1,\cdots,a_\ell,b_\ell,c)\in G^{2\ell+1},
\prod_{i=1}^\ell[a_i,b_i]=c^2 r\right\}
\end{eqnarray*}or

Case 2: $\Si=\Si_\ell\sharp ~\mbox{Klein bottle}$
\begin{eqnarray*}
N_r&=& \left\{(V,c,V,cr)=(a_1,b_1,\cdots,a_\ell,b_\ell,d,c,a_1,b_1,\cdots,a_\ell,b_\ell,d,cr)\mid\right.\\
&& \left. (a_1,b_1,\cdots,a_\ell,b_\ell,d,c)\in
G^{2\ell+2},\prod_{i=1}^\ell[a_i,b_i]=cdc^{-1}dr \right\}
\end{eqnarray*}

In this way, we choose a nice representative $\tilde{x} \in M$ for
any generic point $[x] \in (M/K)^\tau$ where ``generic'' was
defined in Definition \ref{gen}.


Let $P:M\ra M/K$ be the natural projection. We have
\[ S:=\bigcup_{r\in Z(G)} P(N_r)\subset (M/K)^{\tau}, \]
and $S$ is a dense subset of $(M/K)^{\tau}$. Clearly $\bigcup_{r
\in Z(G)}P(N_r) \subset (M/K)^{\tau}$. The reason for $S$ being
dense in $(M/K)^{\tau}$ is that, for generic $[x] \in
(M/K)^{\tau}$, we can always find some $h$ such that $h.x \in N_r$
for some $r \in Z(G)$ (from the previous argument) such that
$P(h.x)=[h.x]=[\tilde{x}]=[x]$, which means $\bigcup_{r \in Z(G)}
P(N_r)$ covers generic points in $(M/K)^{\tau}$.

Thus, there are at most $Z(G)$ copies of $P(N_r)$ which cover
generic points in $(M/K)^{\tau}$. We would like to know if there
is any relation between $P(N_{r_1})$ and $P(N_{r_2})$ for $r_1
\neq r_2$ in $Z(G)$.

\begin{lm}
Let $2Z(G) \subset Z(G)$ denote the subgroup $\{g^2\mid g\in
Z(G)\}$ of the finite abelian group $Z(G)$. If $[r_2]=[r_1] \in
Z(G)/2Z(G)$, then $P(N_{r_1})=P(N_{r_2})$. Moreover, if $[r]\neq[e]$,
then for generic $[x]\in P(N_r)$, $[x]\notin P(N_e)$.
\end{lm}

\paragraph{Proof}

If $[r_2]=[r_1] \in Z(G)/2Z(G)$, then $r_2=r_1 s^2$ for some $s\in Z(G)$, and
\[
(s,e)\cdot(V,c,V,cr_2)=(V,sc,V,scr_1),
\]
which means \[ [(V,c,V,cr_2)]=[(V,sc,V,scr_1)] ,
\]
where the left hand side is a point in $P(N_{r_2})$ and the right
hand side is a point in $P(N_{r_1})$,
 i.e. $P(N_{r_2}) \subset P(N_{r_1})$ which implies
$P(N_{r_1})=P(N_{r_2})$. Thus, if $[r_1]=[r_2] \in Z(G)/2Z(G)$,
then $P(N_{r_1})=P(N_{r_2})$ (in $(M/K)^{\tau}$).

Now let $r$ be an element in $Z(G)$ such that $r\notin 2Z(G)$. We
want to show that \(P(N_r) \cap P(N_e)\) is generically empty,
i.e. for generic $[x]\in P(N_r)$, $[x]\notin P(N_e)$. Let
\[x=(V,c,V,cr)\in N_r.\]

Suppose that $[x]\in P(N_e)$. Then there exists $k=(g_1,g_2)\in K$
such that $k\cdot x\in N_e$. We have

\[
(g_1Vg_1^{-1},g_1cg_2^{-1})=(g_2Vg_2^{-1},g_2cg_1^{-1})\]
$\Longleftrightarrow$
\[ (V,c)
=(g_1^{-1}g_2Vg_2^{-1}g_1, g_1^{-1}g_2 c g_1^{-1}g_2 r)
\]
For generic $[x]\in P(N_r)$ we have $g_1^{-1}g_2=s\in Z(G)$, and
$c=g_1^{-1}g_2 c g_1^{-1}g_2 r$ gives us $s^2 r=e$, which
contradicts the fact that $r\notin 2Z(G)$.

$\Box$

This simplifies the relation between $\bigcup_{r \in Z(G)}P(N_r)$
and $(M/K)^{\tau}$ as follows:

\begin{lm}
Given $[r]\in Z(G)/2Z(G)$, set $P_{[r]}=P(N_r)$. We have
\[
I(M^\tau/K^\tau)=P_{[e]}\subset S:=\bigcup_{[r]\in Z(G)/2Z(G)}
P_{[r]}= (M/K)^\tau,
\]
where $e$ is the identity element of $Z(G)$ and the union is
(generically) disjoint union.
\end{lm}

\paragraph{Proof.}
We know that $S$ is a dense subset of $(M/K)^\tau$. The set $S$ is
compact since it is the image of the compact subset $\cup_{r\in
Z(G)} N_r$ of $M$ under the continuous map $p:M\ra M/K$. In
particular, $S$ is closed because $\cup_{r\in Z(G)} N_r$ is. Thus
we have
\[
(M/K)^\tau=S=\bigcup_{[r]\in Z(G)/2Z(G)} P_{[r]},
\]
and $P_{[e]}=I(M^\tau/K^\tau)$

$\Box$

Thus we have the following theorem:
\begin{tm}
The map $I: M^{\tau}/K^{\tau} \ra (M/K)^{\tau}$ is surjective if
and only if $\mid Z(G)\mid$ is odd.
\end{tm}
\paragraph{Proof}
The above Lemma says $(M/K)^\tau=\bigcup_{[r]\in
Z(G)/2Z(G)}P_{[r]}$ where $P_{[e]}$ is the image of
$M^{\tau}/K^{\tau}$ under $I$. (This applies to all points, not
only generic points.). This tells us that when $2Z(G)=Z(G)$ the
map $I$ is surjective (not only that the image is a dense set). On
the other hand, each set $P_{[r]}$ is compact and closed, and they
are generically disjoint and thus the set $P_{[e]}$ need not be
dense in $M^{\tau}/K^{\tau}$ if $Z(G)/2Z(G) \neq \{e\}$. This
implies that the map $I$ is surjective if and only if
$2Z(G)=Z(G)$. This is the same as $\mid Z(G)\mid$ being odd
because $Z(G)$ is a finite abelian group, and since $2Z(G)$ is a
subgroup of $Z(G)$, the order of $2Z(G)$ divides the order of the
group $Z(G)$. This means that $2Z(G)=Z(G)$ if and only if the
order of $\Z_{n_i}$ is never even.

$\Box$

\subsection{Injectivity}

Suppose that we have $x,\bar{x}\in M^\tau$ such that
$[x]=[\bar{x}]\in (M/K)^{\tau}$. Then we have $\bar{x}=k\cdot x$
for some $k=(g_1,g_2)\in K$.

Suppose that
\begin{eqnarray*}
x&=&(V,c,V,c)\\
\bar{x}&=&(\bar{V},\bar{c},\bar{V},\bar{c})
\end{eqnarray*}
We have
\begin{eqnarray*}
(\bar{V},\bar{c})&=&(g_1Vg_1^{-1},g_1cg_2^{-1})\\
&=&(g_2Vg_2^{-1},g_2cg_1^{-1})
\end{eqnarray*}
$\Longleftrightarrow$
\[
(V,c)=(g_1^{-1}g_2Vg_2^{-1}g_1,g_1^{-1}g_2cg_1^{-1}g_2)
\]


For the generic case, we have $g_1^{-1}g_2=r\in Z(G)$ and $r^2=e$(
from $c=g_1^{-1}g_2cg_1^{-1}g_2$), where $e$ is the identity
element of $Z(G)$. If we consider the homomorphism $\phi:Z(G)\ra
Z(G)$ given by $g\mapsto g^2$, then $r\in \mathrm{Ker}\phi$. The
element $k=(g_1,g_2)$ is in $K^\tau$ if and only if $r=e$. Let
$d=|\mathrm{Ker}\phi|=|Z(G)/2Z(G)|$. Then $I$ is generically $d$
to $1$. $I$ is generically injective if and only if $2Z(G)=Z(G)$.
In other words, if $x$ and $\tilde{x}$ in $M^\tau$ maps to the
same class in $(M/K)^\tau$ then $r \in Z(G)$ and $r^2=e$, and the
condition for $[x]=[\tilde{x}]$ in $M^{\tau}/K^{\tau}$ is that
$r=e$. Thus for the same class in $(M/K)^\tau$, we could have
(number of preimages of $e$ in $Z(G)$) classes in $M^\tau/K^\tau$,
i.e. $|Z(G)/2Z(G)|$ classes. We have the following theorem:

\begin{tm}
The map $I:M^\tau/K^\tau \ra (M/K)^\tau$ is generically
$|Z(G)/2Z(G)|$ to $1$. In particular, $I$ is injective at generic
points of $M^{\tau}/K^{\tau}$ if and only if $\mid Z(G)\mid$ is
odd.
\end{tm}
\begin{rem} \end{rem}
As we see now, the topological type of the nonorientable surface
does not affect the injectivity and surjectivity of the map $I$.

\section{Examples}

Let $I:M^\tau/K^\tau\ra (M/K)^\tau$ be defined as in
Section~\ref{projective}, where $M^\tau/K^\tau$ is the moduli
space of gauge equivalence classes of flat $G$-connections on a
nonorientable compact surface $\Si$, and $(M/K)^\tau$ is the fixed
point set of the involution on the moduli space of gauge
equivalence classes of flat $G$-connections on its orientable
double cover $\tS$. Also, recall that the generic points are the
smooth part of the manifold $M$.

\begin{enumerate}
\item $G=SU(n)$: If $n$ is odd, we have $Z(G)=\Z_n=2Z(G)$, so
      $I$ is surjective and generically injective. If $n$ is even, we have
      $Z(G)/2Z(G)=\{[e],[-e]\}\cong\Z/2\Z$, so $I$ is not surjective and
      is generically $2$ to $1$.
\item $G=Sp(n)$: We have $Z(G)=\{e,-e\}$ and $2Z(G)=\{e\}$, so $I$
is not
      surjective and is generically $2$ to $1$.
\item $G=Spin(n)$: If $n$ is odd, we have $Z(G)=\{1,-1\}$ and
$2Z(G)=\{1\}$,
      so $I$ is not surjective and is generically $2$ to $1$.
      If $n$ is even, we have
      $Z(G)=\{1,-1,e_1e_2\cdots e_n, -e_1e_2\cdots e_n\}
      \cong \Z/2\Z\times \Z/2\Z$, and $2Z(G)=\{1\}$,
      so $I$ is not surjective and is generically $4$ to $1$.
\item $G=SO(n)$: If $n$ is odd, we have $Z(G)=\{e\}=2Z(G)$, so
      $I$ is surjective and generically injective.
       If $n$ is even, we have $Z(G)=\{e,-e\}$ and
      $2Z(G)=\{e\}$, so $I$ is not surjective and is generically $2$ to $1$.
\end{enumerate}

\section{moduli space of semistable vector bundles on a Riemann surface}

There is a very important relation between the moduli space of
flat $SU(k)$-connections over a Riemann surface and the moduli
space of semistable holomorphic vector bundles of rank $k$, degree
$0$ on a Riemann surface as shown in \cite{NS}. Using our result,
we find an interesting Lagrangian submanifold of the moduli space
of semistable holomorphic vector bundles of rank $k$ and degree
$d$ on a Riemann surface of genus $\ell$,
\[\M(k,d)=
\{(a_1,b_1,\cdots,a_\ell,b_\ell)\in G^{2\ell}\mid
~\prod_{i=1}^{\ell} [a_i,b_i]=\exp^{2\pi id/k}\mathbb{I}_k \}/G,
\]

For $k=2n+1$, the moduli space of flat $SU(2n+1)$-connections on a
non-orientable surface is a minimal Lagrangian submanifold of the
moduli space of flat $SU(2n+1)$-connections on its orientable
double cover which is a Riemann surface. Thus the moduli space of
flat $SU(2n+1)$-connections on a nonorientable surface is a
minimal Lagrangian submanifold of the moduli space of semistable
holomorphic vector bundles of rank $2n+1$ and degree $0$ on its
orientable double cover.

For $k=2n$, we know that the moduli space of flat
$PSU(2n)$-connections on a non-orientable surface is a minimal
Lagrangian submanifold of the moduli space of flat
$PSU(2n)$-connections on its orientable double cover because the
center of $PSU(2n)$ is the trivial group. We look at the
commutative diagram below,
\[
\begin{CD}
SU(2n)^{2\ell} @>{\tilde{\mu}}>> SU(2n)\\
@V{\pi^{2\ell}}VV @VV{\pi}V\\
PSU(2n)^{2\ell} @>{\mu}>> PSU(2n)
\end{CD}
\]where $\pi$ is the projection from $SU(2n)$ to $PSU(2n)$ and $\mu$ is the
product of commutators. If $h\in PSU(2n)^{2\ell}$ is an element
such that $\mu(h)=e$ the identity element in $PSU(2n)$, then there
is a lifting $\tilde{h}\in SU(2n)^{2\ell}$ of $h$ such that the
whole diagram commutes. Thus, $\tilde{\mu}(\tilde{h})\in
Ker\pi=Z(SU(2n))$. Thus $\tilde{\mu}(\tilde{h})=\exp^{2\pi
id/2n}\mathbb{I}_{2n}$ for some integer $0\leq d\leq 2n$ which
means $\tilde{h}$ represents a semistable vector bundle of rank
$2n$ and degree $d$. The question is which $(d,2n)$ this
$\tilde{h}$ represents.

First, we need to go back to two earlier papers \cite{Li}\cite{hl}
about the connected components of the moduli spaces. In \cite{Li}
we see that the moduli space of flat $PSU(2n)$-connections on a
Riemann surface has $2n$ connected components since the order of
$\pi_1(PSU(2n))$ is $2n$. Thus the moduli space of flat
$PSU(2n)$-connections on a Riemann surface (equivalently
representations of $\pi_1(\Si)$ into $PSU(2n)$ mod conjugacy) is
equivalent under $\pi^{2\ell}$ to the disjoint union of the moduli
spaces $\M(2n,d)$ of semistable holomorphic vector bundles of rank
$2n$ and of degree $d$, $d$ from $0$ to $2n-1$, by the
Narasimhan-Seshadri theorem \cite{NS}. In \cite{hl}, we can see
that the moduli space of flat $PSU(2n)$-connections on a
nonorientable surface has only two connected components since the
order of $\pi_1(PSU(2n))/2\pi_1(PSU(2n))$ is $2$.

Let us recall the map from our model of two base points for the
moduli space to the standard model of one base point for the
moduli space:

Map $\Phi$ for genus $2\ell$,
\[
\Phi:
\{(a_1,b_1,\cdots,a_\ell,b_\ell,c,a'_1,b'_1,\cdots,a'_\ell,b'_\ell,c')\in
G^{4\ell+2}\mid \prod_{i=1}^{\ell}[a_i,b_i]=cc',
\prod_{i=1}^{\ell}[a'_i,b'_i]=c'c\}/G\times G
 \]\[\ra \{(a_1,b_1,\cdots,a_\ell,b_\ell,a'_1,b'_1,\cdots,a'_\ell,b'_\ell)\in G^{4\ell} \mid \prod_{i=1}^{\ell}[a_i,b_i]=e\}/G
\]
is defined as
$$\Phi:
~~~[(a_1,b_1,\cdots,a_\ell,b_\ell,c,a'_1,b'_1,\cdots,a'_\ell,b'_\ell,c')]~~~~~~~~$$
\begin{eqnarray}\label{phimap}
&&~~~~~~~~~~~\mapsto [(a_1,b_1,\cdots,a_\ell,b_\ell,c'^{-1}b'_\ell
c',c'^{-1}a'_\ell c',\cdots,c'^{-1}b'_1c',c'^{-1}a'_1c')].
\end{eqnarray}

Notice that this map is a well defined bijection:

If\begin{eqnarray*}\lefteqn{
(a_1,b_1,\cdots,a_\ell,b_\ell,c,a'_1,b'_1,\cdots,a'_\ell,b'_\ell,c')}\\&\stackrel{\Phi}{\mapsto}&
(a_1,b_1,\cdots,a_\ell,b_\ell,c'^{-1}b'_\ell c',c'^{-1}a'_\ell
c',\cdots,c'^{-1}b'_1c',c'^{-1}a'_1c'),\end{eqnarray*} then
\begin{eqnarray*}
\lefteqn{
(e,g).(a_1,b_1,\cdots,a_\ell,b_\ell,c,a'_1,b'_1,\cdots,a'_\ell,b'_\ell,c')}\\
&=&(a_1,b_1,\cdots,a_\ell,b_\ell,cg^{-1},ga'_1g^{-1},gb'_1g^{-1},\cdots,ga'_\ell g^{-1},gb'_\ell g^{-1},gc')\\
 &\stackrel{\Phi}{\mapsto}&
(a_1,b_1,\cdots,a_\ell,b_\ell,(gc')^{-1}gb'_\ell
g^{-1}(gc'),(gc')^{-1}ga'_\ell g^{-1}(gc'),\\&&
\cdots,(gc')^{-1}gb'_1 g^{-1}(gc'),(gc')^{-1}ga'_1g^{-1}(gc'))\\
&=&(a_1,b_1,\cdots,a_\ell,b_\ell,c'^{-1}b'_\ell c',c'^{-1}a'_\ell
c',\cdots,c'^{-1}b'_1c',c'^{-1}a'_1c')\end{eqnarray*} which are
the same elements, and \begin{eqnarray*}\lefteqn{
(g,e).(a_1,b_1,\cdots,a_\ell,b_\ell,c,a'_1,b'_1,\cdots,a'_\ell,b'_\ell,c')}\\
&=&(ga_1g^{-1},gb_1g^{-1},\cdots,ga_\ell g^{-1},gb_\ell g^{-1},gc,a'_1,b'_1,\cdots,a_\ell,b_\ell,c'g^{-1})\\
&\stackrel{\Phi}{\mapsto}& (ga_1g^{-1},gb_1g^{-1},\cdots,ga_\ell
g^{-1},gb_\ell g^{-1},(c'g^{-1})^{-1}b'_\ell(c'g^{-1}),\\&&
(c'g^{-1})^{-1}a'_\ell(c'g^{-1}),\cdots,(c'g^{-1})^{-1}b'_1(c'g^{-1}),
(c'g^{-1})^{-1}a'_1(c'g^{-1}))\\
&=&g(a_1,b_1,\cdots,a_\ell,b_\ell,c'^{-1}b'_\ell c',c'^{-1}a'_\ell
c',\cdots,c'^{-1}b'_1c',c'^{-1}a'_1c')g^{-1}\end{eqnarray*} which
means they are members of the same equivalence classes in the
moduli space. Thus
\begin{eqnarray*}\lefteqn{
[(a_1,b_1,\cdots,a_\ell,b_\ell,c,a'_1,b'_1,\cdots,a'_\ell,b'_\ell,c')]
}\\&\stackrel{\Phi}{\mapsto}&
[(a_1,b_1,\cdots,a_\ell,b_\ell,c'^{-1}b'_\ell c',c'^{-1}a'_\ell
c',\cdots,c'^{-1}b'_1c',c'^{-1}a'_1c')].
\end{eqnarray*}is well defined.

According to the argument appeared in \cite{h}, the moduli space
$\M_{PSU(2n)}$ of flat $PSU(2n)$-connections on the connected sum
of $\Si_\ell$ and $RP^2$ has two connected components:
\begin{eqnarray*}
\M_{PSU(2n)}&=& \{(a_1,b_1,\cdots,a_\ell,b_\ell,c)\in
PSU(2n)^{2\ell+1}\mid
\prod_{i=1}^{\ell}[a_i,b_i]c^2=e\}/PSU(2n)\\
 &=&\pi^{2\ell+1}(\{(\tilde{a_1},\tilde{b_1},\cdots,\tilde{a_\ell},\tilde{b_\ell},\tilde{c})
 \in SU(2n)^{2\ell+1}\mid
\prod_{i=1}^{\ell}[\tilde{a_i},\tilde{b_i}]\tilde{c}^2=e\}/SU(2n))\\
&\bigcup&
\pi^{2\ell+1}(\{(\tilde{a_1},\tilde{b_1},\cdots,\tilde{a_\ell},\tilde{b_\ell},\tilde{c})
 \in SU(2n)^{2\ell+1}\mid
\prod_{i=1}^{\ell}[\tilde{a_i},\tilde{b_i}]\tilde{c}^2=\exp^{\frac{2\pi
i}{2n}}\mathbb{I}_{2n} \}/SU(2n))
\end{eqnarray*}
 where the map
\[\pi: SU(2n) \mapsto PSU(2n)=SU(2n)/Z(SU(2n))\] is the projection
map.

If we look at this moduli space $\M_{PSU(2n)}$ on the connected
sum of $\Si_\ell$ and $RP$ as the fixed point set of the moduli
space $\tilde{\M}_{PSU(2n)}$ of flat $PSU(2n)$ connections on
$\Si_{2\ell}$, it becomes
\begin{eqnarray*}
\M_{PSU(2n)}&=&\{(a_1,b_1,\cdots,a_\ell,b_\ell,c,a_1,b_1,\cdots,a_\ell,b_\ell,c)\in
PSU(2n)^{4\ell+2} \mid\\&& \prod_{i=1}^{\ell}[a_i,b_i]c^2=e\}
/(PSU(2n)\times PSU(2n))\\& \subset&
\tilde{\M}_{PSU(2n)}\end{eqnarray*}where
\begin{eqnarray*}
\tilde{\M}_{PSU(2n)}&=&\{(a_1,b_1,\cdots,a_\ell,b_\ell,c,a'_1,b'_1,\cdots,a'_\ell,b'_\ell,c')
\in PSU(2n)^{4\ell+2} \mid\\&&
\prod_{i=1}^{\ell}[a_i,b_i]=cc',\prod_{i=1}^{\ell}[a'_i,b'_i]=c'c\}
/(PSU(2n)\times PSU(2n))
\end{eqnarray*}

So if a point
$(a_1,b_1,\cdots,a_\ell,b_\ell,c,a_1,b_1,\cdots,a_\ell,b_\ell,c)\in
\tilde{\M}_{PSU(2n)}$ is in $\M_{PSU(2n)}$ then it must either
satisfy $\prod_{i=1}^{\ell}[\tilde{a}_i,\tilde{b}_i]\tilde{c}^2=e$
or $\prod_{i=1}^{\ell}[\tilde{a}_i,\tilde{b}_i]\tilde{c}^2=k$ for
some $k \in Z(SU(2n))$ such that $[k] \neq [e] \in
Z(SU(2n))/2Z(SU(2n))$ where $\tilde{a}_i$, $\tilde{b}_i$,
$\tilde{c}$ are the lifting of $a$, $b$, $c$ from $PSU(2n)$ to
$SU(2n)$.

On the other hand, $\tilde{\M}_{PSU(2n)}$ has another description

$$\{(a_1,b_1,\cdots,a_\ell,b_\ell,a'_1,b'_1,\cdots,a'_\ell,b'_\ell) \in
PSU(2n)^{4\ell}\mid
\prod_{i=1}^{\ell}[a_i,b_i]\prod_{i=1}^{\ell}[a'_i,b'_i]=e\}/PSU(2n).$$

which is the standard model for the moduli space with one base
point and these two models are identified by the map $\Phi$ we
defined by equation (\ref{phimap}).

Thus we only need to know what is the product of the commutators
of the lifting in $\tilde{\M}_{SU(2n)}$ of a point in
$\M_{PSU(2n)} \subset \tilde{\M}_{PSU(2n)}$.


Suppose
$(a_1,b_1,\cdots,a_\ell,b_\ell,c,a_1,b_1,\cdots,a_\ell,b_\ell,c)\in
\M_{PSU(2n)} \subset PSU(2n)^{4\ell+2}$. If $a_i$ lifts to
$\tilde{a}_i$, $b_i$ lifts to $\tilde{b}_i$, and $c$ lifts to
$\tilde{c}$, then $
\prod_{i=1}^{\ell}[\tilde{a}_i,\tilde{b}_i]\tilde{c}^2=e $ or $k$
for some $k\in Z(SU(2n))$ such that $[k] \neq [e] \in
Z(SU(2n))/2Z(SU(2n))$ and the map $\Phi$ sends
\begin{eqnarray*} \lefteqn{
(a_1,b_1,\cdots,a_\ell,b_\ell,c,a_1,b_1,\cdots,a_\ell,b_\ell,c)\in
\M_{PSU(2n)} \subset \tilde{M}_{PSU(2n)} } \\ &\mapsto&
(a_1,b_1,\cdots,a_\ell,b_\ell,c^{-1}b_\ell c,c^{-1}a_\ell
c,\cdots,c^{-1}b_1c,c^{-1}a_1c)\in
\tilde{\M}_{PSU(2n)}.\end{eqnarray*}

Thus the product of the commutators of the lifting is \[
\prod_{i=1}^{\ell}[\tilde{a}_i,\tilde{b}_i]\prod_{i=1}^{\ell}[\tilde{c}^{-1}\tilde{b}_i\tilde{c},
\tilde{c}^{-1}\tilde{a}_i\tilde{c}]=\prod_{i=1}^{\ell}[\tilde{a}_i,\tilde{b}_i]\tilde{c}^{-1}
\prod_{i=1}^{\ell}[\tilde{b}_i,\tilde{a}_i]\tilde{c}=\left\{\begin{array}{ll}e&\mbox{if
$\prod_{i=1}^{\ell}[\tilde{a}_i,\tilde{b}_i]\tilde{c}^2=e$
}\\kk^{-1}=e&\mbox{if
$\prod_{i=1}^{\ell}[\tilde{a}_i,\tilde{b}_i]\tilde{c}^2=k$}\end{array}
\right. \]

This means that the lifting of the moduli space of flat
$PSU(2n)$-connections over the connected sum of $\Si_\ell$ and
$RP^2$ is (in) the moduli space of semi-stable vector bundle of
rank $2n$ and degree 0 over $\Si_{2\ell}$. Thus a natural
Lagrangian submanifold of $\M(2n,0)$ over $\Si_{2\ell}$ is
(equivalent to) the moduli space of flat $PSU(2n)$-connections
over the connected sum of $\Si_{2\ell}$ and $RP^2$ and a natural
Lagrangian submanifold of $\M(2n,0)$ over $\Si_{2\ell+1}$ is the
moduli space of flat $PSU(2n)$-connections over the connected sum
of $\Si_{\ell}$ and the $Klein ~bottle$.

Another way to look at the problem is that:

For SU(2n), the fixed point set of the moduli space
$\tilde{\M}_{SU(2n)}$ is isomorphic to the generically disjoint
union of $P_0$ and $P_1$ where
$$P_0=\{(a_1,b_1,\cdots,a_\ell,b_\ell,c,a_1,b_1,\cdots,a_\ell,b_\ell,c)\in SU(2n)^{4\ell+2}\mid
\prod_{i=1}^{\ell}[a_i,b_i]=c^2 \}/SU(2n)\times SU(2n)$$ and
choosing any $r \in Z(SU(2n))$, $[r] \neq [e]$,
$$P_1=\{ (a_1,b_1,\cdots,a_\ell,b_\ell,c,a_1,b_1,\cdots,a_\ell,b_\ell,cr)\in SU(2n)^{4\ell+2}
\mid \prod_{i=1}^{\ell}[a_i,b_i]=c^2r \}/SU(2n)\times SU(2n).$$

The piece $P_0$ is in fact the moduli space of $SU(2n)$-bundles
over the connected sum of $\Si_\ell$ and $RP^2$ which is
$$\{(a_1,b_1,\cdots,a_\ell,b_\ell,c)\in SU(2n)^{2\ell+1}\mid \prod_{i=1}^{\ell}[a_i,b_i]=c^2\}/SU(2n),$$ thus we have
left only $P_1$ unidentified.

If we look at the moduli space of $PSU(2n)$-bundles over the
connected sum of $\Si_\ell$ and $RP^2$, according to \cite{h}, we
also have two pieces
$$\pi^{2\ell+1}(\{(a_1,b_1,\cdots,a_\ell,b_\ell,c)\in SU(2n)^{2\ell+1}
\mid \prod_{i=1}^{\ell}[a_i,b_i]=c^2\}/SU(2n))$$ and
$$\pi^{2\ell+1}(\{(a_1,b_1,\cdots,a_\ell,b_\ell,c)\in SU(2n)^{2\ell+1}
\mid \prod_{i=1}^{\ell}[a_i,b_i]=c^2r\}/SU(2n))$$ for $r \in
Z(SU(2n))$ and $[r] \neq [e] \in Z(SU(2n))/2Z(SU(2n))$, where
$\pi$ is the projection map from $SU(2n)$ to $PSU(2n)$.

Now the first piece is identified with $\pi^{2\ell+1}(P_0)$. What
about the second piece? Note that the space
$$ \{(a_1,b_1,\cdots,a_\ell,b_\ell,c)\in G^{2\ell+1}\mid
\prod_{i=1}^{\ell}[a_i,b_i]=c^2r\}/G$$ is identified with the
space
$$ \{(a_1,b_1,\cdots,a_\ell,b_\ell,c,a_1,b_1,\cdots,a_\ell,b_\ell,cr)\mid
\prod_{i=1}^{\ell}[a_i,b_i]=c^2r\}/G\times G$$ by the isomorphism
$$ (a_1,b_1,\cdots,a_\ell,b_\ell,c) \mapsto
 (a_1,b_1,\cdots,a_\ell,b_\ell,c,a_1,b_1,\cdots,a_\ell,b_\ell,cr)$$
so the second piece can be identified with $\pi^{2\ell+1}(P_1)$.

We can generalize this idea from $G=SU(n)$ to any compact
connected semisimple Lie group $G$, and draw the following
relations:
$$\M_G(\Si_\ell\sharp RP^2)=\{(a_1,b_1,\cdots,a_\ell,b_\ell,c,a_1,b_1,\cdots,a_\ell,b_\ell,c)
\in G^{4\ell+2}\mid\prod_{i=1}^{\ell}[a_i,b_i]=c^2\}/(G\times G)$$
$$\downarrow$$
$$\mbox{the fixed point set of}
~\tilde{\M}_G(\Si_{2\ell})~=~~ P_0 \cup P_1(~\mbox{or more
pieces})$$ where
\begin{eqnarray*}
P_0&=&\{(a_1,b_1,\cdots,a_\ell,b_\ell,c,a_1,b_1,\cdots,a_\ell,b_\ell,c)
\in G^{4\ell+2}\mid \prod_{i=1}^{\ell}[a_i,b_i]=c^2\}/(G\times
G),\\
P_1&=&\{(a_1,b_1,\cdots,a_\ell,b_\ell,c,a_1,b_1,\cdots,a_\ell,b_\ell,cr)
\in G^{4\ell+2}\mid \prod_{i=1}^{\ell}[a_i,b_i]=c^2r\}/(G\times
G),\end{eqnarray*} for $r \in Z(G)$ such that $[r] \neq [e] \in
Z(G)/2Z(G)$, and if we have more pieces then we have $r_i \in
Z(G)$ such that $[r_i] \neq [r_j] \neq [e] \in Z(G)/2Z(G)$.

Since the center of $G/Z(G)$ is trivial, the fixed point set of
$\M_{G/Z(G)}(\Si_{2\ell})$ is equivalent to
$\M_{G/Z(G)}(\Si_\ell\sharp RP^2)$, we have
\begin{eqnarray*}&&\mbox{The fixed point set of}~~\M_{G/Z(G)}(\Si_{2\ell})\\
&=&\{(a_1,b_1,\cdots,a_\ell,b_\ell,c,a_1,b_1,\cdots,a_\ell,b_\ell,c)\in
(G/Z(G))^{4\ell+2}\mid \prod_{i=1}^{\ell}[a_i,b_i]=c^2\}/(G/Z(G)
\times G/Z(G))
\\&=&\{(a_1,b_1,\cdots,a_\ell,b_\ell,c)\in (G/Z(G))^{2\ell+1}\mid
\prod_{i=1}^{\ell}[a_i,b_i]=c^2\}/(G/Z(G))
\\&=&\pi^{2\ell+1}(\{(\tilde{a}_1,\tilde{b}_1,\cdots,\tilde{a}_\ell,\tilde{b}_\ell,\tilde{c})
\in G^{2\ell+1}\mid \prod_{i=1}^{\ell}
[\tilde{a}_i,\tilde{b}_i]=\tilde{c}^2\}/G)\\&&
\cup\pi^{2\ell+1}(\{(\tilde{a}_1,\tilde{b}_1,\cdots,\tilde{a}_\ell,\tilde{b}_\ell,\tilde{c})\in
G^{2\ell+1}\mid
\prod_{i=1}^{\ell}[\tilde{a}_i,\tilde{b}_i]=\tilde{c}^2r\}/G)
\\&=& \pi^{2\ell+1}(P_0) \cup \pi^{2\ell+1}(P_1)
\end{eqnarray*}

{\bf Conclusion:} \begin{enumerate}

\item  For $G$ for which the order of its center is not odd, the
fixed point set of the involution on the moduli space of flat
$G$-bundles over $\Si_{2\ell}$ has more than one piece say
$P_1$(or $P_1,P_2,P_3$) where $P_0$ is identified with the moduli
space of flat $G$-bundles over the connected sum of $\Si_\ell$ and
$RP^2$, then this extra piece $P_1$ or these extra pieces
$P_1,P_2,P_3$ can be identified with the other connected
components of the moduli space of flat $G/Z(G)$-bundles over the
connected sum of $\Si$ and $RP^2$, while the connected component
$\pi^{2\ell+1}(\{(a_1,b_1,\cdots,a_\ell,b_\ell,c)\in
G^{2\ell+1}\mid \prod_{i=1}^{\ell}[a_i,b_i]=c^2\}/G)$ is
identified with $\pi^{2\ell+1}(P_0)$.

or\item  The lifting from $G/Z(G)$ to $G$ of the moduli space of
$G/Z(G)$-bundles over $\Si\sharp RP^2$, $\M_{G/Z(G)}(\Si \sharp
RP^2)$, is the fixed point set of the moduli space of $G$-bundles
over $\tilde{\Si}$, $\M_G(\tilde{\Si})$. In other words, the fixed
point set of the moduli space of $G$-bundles over $\tilde{\Si}$ is
identified (through the projection $\pi$) with the fixed point set
of the moduli space of $G/Z(G)$-bundles over $\tilde{\Si}$.

\end{enumerate}

\begin{rem}\label{CP3}
\end{rem}
We know that the moduli space of flat $SU(2)$-connections on the
genus 2 Riemann surface is the smooth $CP^3$ as showed in
\cite{NR}. We want to show that , $N$, the fixed point set of the
involution which we referred before is the natural $RP^3$ with the
metric on $RP^3$ inherited from the Fubini-Study metric on $CP^3$.

Notice that $RP^3$ is a totally real, totally geodesic submanifold
in $CP^3$. Fix one point $p \in CP^3$, the tangent space at $p$
decomposes as $T_pCP^3=T_pRP^3 \oplus JT_pRP^3$. Since $N$ and
$RP^3$ are both real in $CP^3$, it is possible to send $T_qN$ to
$T_pRP^3$. In fact, we can find an element in $PSU(4)$ the
isometry group of $CP^3$ such that it sends $T_qN$ to $T_pRP^3$.


It suffices to show that any element of $GL(3,\mathbb{C})$ maps
one totally real vector space in $\mathbb{C}^3$ to another.
Suppose $W$ is a totally real vector space in $\mathbb{C}$. Let
$\{w_1,w_2,w_3\}$ denote a basis of $W$. Then $\{Jw_1,Jw_2,Jw_3\}$
is a basis of $JW$ and $W\cap JW=\{0\}$. Suppose $k$ is any
element of $GL(3,\mathbb{C})$. We want to show that $kW\cap
JkW=\{0\}$. Suppose $p \in kW$ then $p=akw_1+bkw_2+ckw_3$ for some
constant $a,b,c$. Thus $p=k(aw_1+bw_2+cw_3)$. Now if $p \in JkW$
then $p=AJkw_1+BJkw_2+CJkw_3$ for some constant $A,B,C$. Recall
that $GL(3,\mathbb{C})=\{k\in GL(6,\mathbb{R}), kJ=Jk\}$. Thus
$p=k(AJw_1+BJw_2+CJw_3)$. Since $\mbox{det}(k)\neq 0$, we have
$aw_1+bw_2+cw_3=AJw_1+BJw_2+CJw_3$ which gives all the constants
$0$ because of $W\cap JW=\{0\}$.

Now apply the exponential map to both $T_pN$ and $T_qRP^3$, by the
uniqueness of geodesics, $N \simeq RP^3$ and the metric on $N$ is
the induced metric on $RP^3$ since we know $\M \simeq CP^3$ with
the standard Fubini-Study metric. The reason that $N$ is smooth is
that $N$ is the same as $RP^3$ by the exponential map, so it is
smooth because $RP^3$ is. $N$ is not a proper subset of $RP^3$
because $N$ is also compact and has the same dimension as $RP^3$
so it is not possible for it to be a proper subset.

\end{document}